\documentclass[11pt, reqno, a4paper]{amsart}
\usepackage{a4wide}
\usepackage{amsmath, amsthm, amssymb}
\usepackage{enumerate}
\usepackage[all]{xy}

\usepackage{color}
\usepackage{ifpdf}
\ifpdf 
  \usepackage[pdftex]{graphicx}
  \DeclareGraphicsExtensions{.pdf,.png,.mps}
  \usepackage{pgf}
  \usepackage{tikz}
\else 
  \usepackage{graphicx}
  \DeclareGraphicsExtensions{.eps,.bmp}
  \usepackage{pgf}
  \usepackage{tikz}
  \usepackage{pstricks}
\fi
\usepackage{epic}
\usepackage{wrapfig}

\newlength{\sperr}

\newtheorem{thm}{Theorem}
\newtheorem{cor}[thm]{Corollary}
\newtheorem{corollary}[thm]{Corollary}

\theoremstyle{definition}
\newtheorem{defn}[thm]{Definition}
\newtheorem*{rmk}{Remark}
\newtheorem*{ex}{Example}

\newcommand\set[1]{\left\{#1\right\}}

\newcommand\floor[1]{\left\lfloor #1 \right\rfloor}

\mathchardef\mhyphen="2D
\def\Z{\mathbb{Z}}

\DeclareMathOperator\fix{fix}
\DeclareMathOperator\exc{exc}
\DeclareMathOperator\wex{wex}
\DeclareMathOperator\maj{maj}
\DeclareMathOperator\inv{inv}
\DeclareMathOperator\Inv{Inv}
\DeclareMathOperator\cros{cros}
\DeclareMathOperator\nest{nest}
\DeclareMathOperator\fmax{fmax}

\DeclareMathOperator\des{des}
\DeclareMathOperator\ndes{ndes}

\DeclareMathOperator\MAD{MAD}
\DeclareMathOperator\adj{adj}
\DeclareMathOperator\suc{suc}

\DeclareMathOperator\toht{(31\mhyphen 2)}
\DeclareMathOperator\thto{(2\mhyphen 31)}
\DeclareMathOperator\thot{(2\mhyphen 13)}
\def\ptoht{31\mhyphen 2}
\def\pthto{2\mhyphen 31}
\def\pthot{2\mhyphen 13}

\def\S{\mathfrak{S}}
\def\D{\mathfrak{D}}
\def\A{\mathfrak{A}}
\def\L{\mathfrak{L}}
\def\M{\mathfrak{M}}
\def\DP{\mathcal{D}}
\def\DPD{\mathcal{P}}

\title{The $q$-tangent and $q$-secant numbers via continued fractions}
\author{Heesung Shin}
\address[Heesung Shin]{Universit\'{e} de Lyon; Universit\'{e} Lyon 1; Institut Camille Jordan; UMR 5208 du CNRS; 43, boulevard du 11 novembre 1918, F-69622 Villeurbanne Cedex, France}
\email{hshin@math.univ-lyon1.fr}
\urladdr{http://math.univ-lyon1.fr/~hshin/}

\author{Jiang Zeng}
\address[Jiang Zeng]{Universit\'{e} de Lyon; Universit\'{e} Lyon 1; Institut Camille Jordan; UMR 5208 du CNRS; 43, boulevard du 11 novembre 1918, F-69622 Villeurbanne Cedex, France}
\email{zeng@math.univ-lyon1.fr}
\urladdr{http://math.univ-lyon1.fr/~zeng/}

\date{\today}



\begin{document}
\maketitle
\begin{abstract}
It is well known that the $(-1)$-evaluation of  the enumerator polynomials of permutations (resp. derangements)
by the number of excedances  gives rise to tangent   numbers (resp. secant numbers).
 Recently, two distinct  $q$-analogues of  the latter result have been discovered by Foata and Han,   and Josuat-Verg\`es, respectively.
  In this paper, we will prove some general  continued fractions expansions formulae,
   which permits us to
  give a unified treatment of Josuat-Verg\`es' two formulae and also to
  derive a new  $q$-analogue of the aforementioned formulae.
 Our approach is based on
a $(p,q)$-analogue of  tangent and secant numbers via  continued fractions  and also
  the generating function of permutations with respect to
the quintuple  statistic
consisting of fixed point number, weak excedance number, crossing number, nesting number and inversion number.
We also give  a combinatorial proof of Josuat-Verg\`es'  formulae by using a new linear model of derangements.
\end{abstract}
\tableofcontents
\section{Introduction}
The {\em tangent numbers $E_{2n+1}$} ($n\geq 0$) and {\em secant numbers $E_{2n}$} ($n\geq 1$)   are the coefficients in the
Taylor expansions of the  tangent  and secant:
\begin{align}\label{eq1}
\tan t&=t+2\frac{t^3}{3!}+16 \frac{t^5}{5!}+272\frac{t^7}{7!}+\cdots +E_{2n+1} \frac{t^{2n+1}}{(2n+1)!} +\cdots,\\
\sec t&=1+\frac{t^2}{2!}+5\frac{t^4}{4!}+61\frac{t^6}{6!}+1385\frac{t^8}{8!}+\cdots+E_{2n} \frac{t^{2n}}{(2n)!} +\cdots.\label{eq2}
\end{align}
Let $\S_n$ be the set of permutations of $[n]:=\{1,2,\ldots, n\}$.  A permutation $\sigma=\sigma_1\sigma_2\cdots \sigma_n\in \S_n$ is
 {\em alternating} (resp. {\em falling alternating}) permutation if
$\sigma_1<\sigma_2$, $\sigma_2>\sigma_3$, $\sigma_3<\sigma_4$, etc. (resp. $\sigma_1>\sigma_2$, $\sigma_2<\sigma_3$, $\sigma_3>\sigma_4$, etc.).
As  proved by  Andr\'e~\cite{And79} in  1879, the number of
alternating permutations in $\S_n$ is $E_n$.
Furthermore,   an integer $i\in [n]$ is called a {\em fixed point}
 (resp.   {\em weak excedance}, {\em excedance}) of $\sigma$  if $ \sigma_i=i$ (resp.  $\sigma_i \ge i$,  $\sigma_i > i$).
Denote the number of fixed points
(resp.  weak excedances, excedances) of $\sigma$  by  $\fix \sigma$ (resp. $\wex\sigma$,  $\exc\sigma$).
The joint distribution of $(\exc, \fix)$ was
first computed by Foata and Sch{\"u}tzenberger~\cite{FS70}, who proved the following formula (see also \cite{KZ04} for another proof):
\begin{equation}\label{eq3}
\begin{aligned}
\sum_{n\geq 0}\frac{t^n}{n!}\sum_{\sigma\in \S_n}x^{\exc \sigma}y^{\fix\sigma}
=&1+yt+(x+{y}^{2})\frac{t^2}{2!}
+({x}^{2}+(3y+1)x+{y}^{3})\frac{t^3}{3!}
+\cdots\\
=&\frac{(1-x)\exp(yt)}{\exp(xt)-x\exp(t)}.
\end{aligned}
\end{equation}
Let $\D_n $  be the  set  of  \emph{derangements} on  $[n]$, i.e.,  permutations without fixed point. From \eqref{eq1}, \eqref{eq2} and \eqref{eq3}
we derive readily    the following combinatorial interpretations of two identities due to Euler and Roselle (see \cite{Eu55,Ros68}):
for $n\ge 1$,
\begin{align}
\label{eq:Euler1}
\sum_{\sigma\in \S_n}(-1)^{\exc \sigma}
&=
\begin{cases}
0 &\text{if $n$ is even},\\
(-1)^{\frac{n-1}{2}}E_n &\text{if $n$ is odd};
\end{cases}\\
\intertext{and}
\label{eq:Euler2}
\sum_{\sigma\in \D_n}(-1)^{\exc \sigma}
&=
\begin{cases}
(-1)^{\frac{n}{2}}E_n &\text{if $n$ is even},\\
0 &\text{if $n$ is odd}.
\end{cases}
\end{align}
The {\em major index}  and  {\em inversion number}
of $\sigma\in \S_n$ are  defined, respectively, by
\begin{align*}
\maj\sigma=\sum_{i=1}^{n-1}i \chi(\sigma_i>\sigma_{i+1}),\quad
\inv\sigma=\sum_{1\leq i<j\leq n} \chi( \sigma_i> \sigma_j),
\end{align*}
where $\chi(A)=1$ if $A$ is true and 0 otherwise.
Let $\A_n^*$ (resp. $\A_n$) be the set of alternating  (resp. {\em falling alternating})  permutations on  $[n]$.
Foata and Han~\cite{FH10} showed  that one can easily derive a $q$-analogue of \eqref{eq:Euler1} and \eqref{eq:Euler2} from
the  generating function for the triple statistic $(\exc, \fix,\maj)$ due to Sherashian and Wachs~\cite{SWachs07}.

\begin{thm}[Foata-Han]
\label{thm:fh} For $n\geq 1$, we have
\begin{align}
\label{eq:qEuler1}
\sum_{\sigma\in \S_n}(-1)^{\exc \sigma} q^{\maj\sigma}
&=
\begin{cases}
0 &\text{if $n$ is even},\\
(-1)^{\frac{n-1}{2}} \sum_{\sigma\in \A_{n}^*} q^{ \inv\sigma} &\text{if $n$ is odd};
\end{cases}\\
\intertext{and}
\label{eq:qEuler2}
\sum_{\sigma\in \D_n}(-1/q)^{\exc \sigma}q^{\maj\sigma}
&=
\begin{cases}
(-1/q)^{\frac{n}{2}} \sum_{\sigma\in \A_{n}^*} q^{ \inv\sigma} &\textrm{if $n$ is even},\\
0 &\textrm{if $n$ is odd}.
\end{cases}
\end{align}
\end{thm}
As noted in \cite{HRZ99,Che08},   $q$-analogues of the tangent and secant numbers
arise also from the following continued fraction expansions~\cite{Fla80}:
\begin{align}
\sum_{n=0}^{\infty} E_{2n+1} t^{2n+1}
&= t+2 t^3 + 16 t^5+ 272 t^7+\cdots
= \cfrac{t}{1-\cfrac{1\cdot 2 t^2}{1-\cfrac{2\cdot 3 t^2}{1-\cfrac{3\cdot 4 t^2}{\ddots}}}},
\label{eq:cf-tangent}
\intertext{and}
\sum_{n=0}^{\infty} E_{2n} t^{2n}
&=1 + t^2+5 t^4+61 t^6 + \cdots = \cfrac{1}{1-\cfrac{ 1^2t^2}{1-\cfrac{2^2t^2}{1-\cfrac{3^2t^2}{\ddots}}}}.
\label{eq:cf-secant}
\end{align}
The  {\em crossing}  number and {\em nesting} number of  $\sigma\in \S_n$ are defined by
\begin{align*}
\cros(\sigma)&=\#\{(i,j):  i < j \le \sigma_i <\sigma_j \quad\text{or}\quad \sigma_i <\sigma_j <i < j\},\\
\nest(\sigma)&=\#\{(i,j):  i < j \le \sigma_j <\sigma_i \quad\text{or}\quad \sigma_j <\sigma_i <i < j \}.
\end{align*}
The number of occurrences of generalized patterns $\ptoht$, $\pthto$, and $\pthot$ in
  $\sigma\in \S_n$  are defined by
\begin{align*}
\toht\sigma&=\#\{(i,j):  i+1<j\quad \text{and}\quad \sigma_i>\sigma_j>\sigma_{i+1}\},\\
\thto\sigma&=\#\{(i,j):  j<i-1\quad \text{and}\quad \sigma_{i-1}>\sigma_j>\sigma_i\},\\
\thot\sigma&=\#\{(i,j): j<i-1\quad \text{and}\quad \sigma_{i-1}<\sigma_k<\sigma_i \}.
\end{align*}
Inspired by the recent $\cros$-analogue of Eulerian polynomials~\cite{Wil05,Cor05} and $\toht$-analogue of
tangent and secant  numbers~\cite{HRZ99,Che08},
 Josuat-Verg\`es~\cite{JV09}  proved
 another   $q$-analogue of
\eqref{eq:Euler1} and \eqref{eq:Euler2}.

\begin{thm}[Josuat-Verg\`es]
\label{thm:JV}
For $n\geq 1$, we have
\begin{align}
\sum_{\sigma\in \S_n}(-1)^{\wex \sigma } q^{\cros\sigma}
&=
\begin{cases}
0&\textrm{if $n$ is even},\\
(-1)^{\frac{n+1}{2}} \sum_{\sigma\in \A_{n}} q^{ \toht\sigma}&\textrm{if $n$ is odd};
\end{cases}
\label{eq:jv1}\\
\intertext{and}
\sum_{\sigma\in \D_n}(-1/q)^{\exc \sigma}q^{\cros\sigma}
&=
\begin{cases}
(-1/q)^{\frac{n}{2}} \sum_{\sigma\in \A_{n}} q^{ \toht\sigma} &\textrm{if $n$ is even},\\
0 &\textrm{if $n$ is odd}.
\end{cases}
\label{eq:jv2}
\end{align}
\end{thm}
\begin{rmk}
Since  $\exc(\sigma) = n-\wex(\sigma^{-1})$ for $\sigma \in \S_n$,
Eq.~\eqref{eq:jv1} reduces to  \eqref{eq:Euler1}  when  $q=1$.
\end{rmk}

Since  \eqref{eq:jv1} and \eqref{eq:jv2} were proved  by expanding   the
generating functions into continued fractions and then identifying the continued fractions in \cite{JV09},
our first aim is then to seek for a
unified proof of Josuat-Verg\`es' theorem from
a continued fraction formula for the generating function of  permutations with respect to
the triple statistic $(\wex,\fix,\cros)$, which should play the role of the Sherashian-Wachs formula in Foata-Han's proof.
Actually we shall prove some more general results (see Theorems~\ref{thm:pq}, \ref{thm:preserving} and \ref{thm:continued_fraction}), which permit us
to obtain  a  new $q$-analogue of \eqref{eq:Euler1} and \eqref{eq:Euler2}.
\begin{thm}
\label{thm:shin-zeng}
For $n\geq 1$, we have
\begin{align}
\sum_{\sigma\in \S_n}(-1/q)^{\exc \sigma } q^{\inv \sigma}
&=
\begin{cases}
0&\textrm{if $n$ is even},\\
(-1)^{\frac{n-1}{2}} \sum_{\sigma\in \A_{n}} q^{ \toht\sigma + 2 \thot\sigma}&\textrm{if $n$ is odd};
\end{cases}
\label{eq:shin-zeng1}\\
\intertext{and}
\sum_{\sigma\in \D_n}(-1)^{\exc \sigma}q^{\inv\sigma}
&=
\begin{cases}
(-q)^{\frac{n}{2}} \sum_{\sigma\in \A_{n}} q^{ \toht\sigma + 2 \thto\sigma} &\textrm{if $n$ is even},\\
0 &\textrm{if $n$ is odd}.
\end{cases}
\label{eq:shin-zeng2}
\end{align}
\end{thm}

Apart from the aforementioned analytic proof, Foata and Han~\cite{FH10} gave also
a combinatorial proof of Theorem~\ref{thm:fh} making use of a linear model of  derangements, called \emph{desarrangements}.
Similarly, we will   provide a combinatorial proof of Theorem~\ref{thm:JV} using  a new linear model of derangements,
i.e., called {\em coderangements}, which are permutations such that every left-to-right
maximum is a descent, see Definition~\ref{def:dz}.

The rest of this paper is organized as follows.
In Section~\ref{sec:2} we define a $(p,q)$-analogue of Tangent and secant numbers as a refined counting of alternating permutations
and show that their generating functions have nice continued fraction expansions. In Sections~\ref{sec:bijection},
by applying  a bijection  similar to that of Clarke-Steingr\'{\i}msson-Zeng~ \cite{CSZ97},
we obtain  another quintuple statistic equidistributed with $(\wex,\fix, \cros, \nest,\inv)$, which results  in another combinatorial interpretation of Theorem~\ref{thm:JV} (see \eqref{eq:sz1} and \eqref{eq:sz2}),
that involves {\em coderangements} (see Definition~\ref{def:dz}).
In Section~\ref{sec:3}
we will prove a continued fraction expansion   for  the generating function
of permutations with respect to
the quadruple statistic $(\wex,\fix, \cros, \nest,\inv)$.
Two  special cases of this formula  will play the role of two $q$-analogues of Roselle's formula.
We then derive an analytic proof of Theorems~\ref{thm:JV} and \ref{thm:shin-zeng} in Section~\ref{sec:analytic}.
 In Section~\ref{sec:combinatorial},
we will construct two killing involutions using the latter interpretations
to prove Theorems~\ref{thm:JV} as well as  a combinatorial proof for \eqref{eq:shin-zeng2}.
Finally, we give a parity independent formula for the $q$-Euler numbers $E_n(q)$ in Section~\ref{sec:parity-indeperndent}.


\section{A $(p,q)$-analogue of tangent and secant numbers}\label{sec:2}
Since the alternating permutations in $\A_n^*$ and the falling alternating permutations in $\A_n$ are  obviously  isomorphic by
$
\sigma_1\sigma_2\cdots \sigma_n\mapsto  (n+1-\sigma_1)(n+1-\sigma_2)\cdots (n+1-\sigma_n)$,
we define the $(p,q)$-analogue of tangent and secant   numbers  by
\begin{align}
E_{2n+1}(p,q) &= \sum_{\sigma \in \A_{2n+1}} p^{\thot\sigma}q^{\toht\sigma},\label{eq:qtang1}
\intertext{and}
E_{2n}(p,q) &= \sum_{\sigma \in \A_{2n}}  p^{\thto\sigma}q^{\toht\sigma}.\label{eq:qsec1}
\end{align}

It turns out that  the generating functions of
$E_n(p,q)$   give rise to remarkable $(p,q)$-analogue of the two continued fractions \eqref{eq:cf-tangent}
and \eqref{eq:cf-secant}.

We first recall some  well-known facts about the combinatorial theory of continued fractions.
A {\em Motzkin path} of length $n$ is a sequence of lattice points $s_i=(i,y_i)$ ($i=0,\ldots,n$) in the plane $\Z^2$ such that
\begin{itemize}
\item $y_i\ge 0$ for all $1\le i\le n-1$ and $y_0=y_n=0$.
\item $y_i-y_{i-1} \in \set{-1,0,1}$.
\end{itemize}
For $i=1,\ldots,n$,  if $y_i-y_{i-1}=-1,0$ and 1 we shall say that  the
 {\em step $(s_{i-1}, s_i)$}  is of type  $\searrow$, $\longrightarrow$, and $\nearrow$, respectively,
and call $y_{i-1}$ the  {\em height} of the step $(s_{i-1}, s_i)$.
Now,
weight  each step $(s_{i-1},s_{i})$ at height $h$ by $w(s_{i-1},s_{i}) = a_h$ (resp. $b_h$, $c_h$) if  the step $(s_{i-1},s_{i})$ is $\nearrow$ (resp. $\longrightarrow$, $\searrow$) and define the weight of the Motzkin path $\gamma=s_0, s_1, \ldots, s_n$ by
$w(\gamma) = \prod_{i=1}^n w(s_{i-1},s_i)$.
A {\em Dyck path} is a Motzkin path without step of type $\longrightarrow$. So the length of a Dyck path must be even.
Let $\M_n$ (resp. $\DP_{2n}$) denote the set of Motzkin  (resp. Dyck)
paths of length $n$ (resp. ${2n}$).
Then,  it is well known (see \cite{Fla80}) that
\begin{align}\label{eq:motzkin}
1+\sum_{n\geq 1}\sum_{\gamma \in \M_n} w(\gamma) t^n
&=\cfrac{1}{1-b_0t-\cfrac{a_0c_1t^2}{1-b_1t-\cfrac{a_1c_2t^2}{\ddots}}},
\intertext{and}
\label{eq:dyck}
1+\sum_{{n}\geq 1}\sum_{\gamma \in \DP_{2n}} w(\gamma) t^{2n}
&=\cfrac{1}{1-\cfrac{a_0c_1t^2}{1-\cfrac{a_1c_2t^2}{1-\cfrac{a_2c_3t^2}{\ddots}}}}.
\end{align}

\begin{thm}\label{thm:pq}
The generating functions for $E_{2n+1}(p,q)$ and $E_{2n}(p,q)$ have the following continued fraction expansions:
\begin{align}
\sum_{n=0}^{\infty} E_{2n+1}(p,q) t^{2n+1}
&= t+ \left( p+q \right) {t}^{3} \notag \\
&\quad + \left( {p}^{4}+3\,q{p}^{3}+4\,{p}^{2}{
q}^{2}+3\,{q}^{3}p+{q}^{4}+{p}^{2}+2\,qp+{q}^{2} \right) {t}^{5} +\cdots \notag\\
&=\cfrac{t}{1-\cfrac{[1]_{p,q}[2]_{p,q}t^2}{1-\cfrac{[2]_{p,q}[3]_{p,q}t^2}{1-\cfrac{[3]_{p,q}[4]_{p,q}t^2}{\ddots}}}},
\label{eq:new-q-tangent2}
\intertext{and}
\sum_{n=0}^{\infty} E_{2n}(p,q) t^{2n}
&=1+{t}^{2}+ \left( {p}^{2}+2\,qp+{q}^{2}+1 \right) {t}^{4} + \cdots \notag\\
&=\cfrac{1}{1-\cfrac{[1]_{p,q}^2t^2}{1-\cfrac{[2]_{p,q}^2t^2}{1-\cfrac{[3]_{p,q}^2t^2}{\ddots}}}}.
\label{eq:new-q-secant2}
\end{align}
where  $[n]_{p,q}=(p^n-q^n)/(p-q)$ for  $n\geq 0$.
\end{thm}
\begin{proof}

A {\em Dyck path diagramme} of length ${2n}$ is a pair  $(\gamma; \xi)$, where $\gamma$ is a Dyck path of length ${2n}$ and $\xi=\xi_1,\ldots,\xi_{2n}$ is a sequence such that $\xi_k \in \{0, 1, \ldots, h_k\}$ if the step $(s_{k-1},s_{k})$ is at height $h_k$. Let $\DPD_{2n}$ be the set of Dyck path diagrammes of length ${2n}$.
If  $(\gamma;\xi)\in \DPD_{2n}$, for  each  step  $(s_{k-1},s_k)$ at height $h_k$ and $\xi_k$ ($1\leq k\leq {2n}$),  we associate the valuation
$$
v((s_{k-1},s_k), \xi_k) = p^{h_k-\xi_k} q^{\xi_k},
$$
and define the valuation of the Dyck path diagramme $(\gamma; \xi)$ by
$\prod_{k=1}^{2n}v((s_{k-1},s_k), \xi_k)$.
Then
\eqref {eq:new-q-tangent2} is equivalent to
\begin{align}\label{eq:FV1}
\sum_{\sigma \in \A_{2n+1}} p^{\thot\sigma}q^{\toht\sigma} = \sum_{(\gamma;\xi) \in \DPD_{2n}}\prod_{k=1}^{2n}v((s_{k-1},s_k), \xi_k)\qquad (n\geq 1).
\end{align}

Similarly, we define the set
$$\DPD^*_{2n}=\set{(\gamma; \xi) \in \DPD_{2n}: (s_{k-1},s_{k})=\nearrow \text{ or } \xi_k < h_k \text{ for all $k=1,\dots,2n$}}.$$
If  $(\gamma;\xi)\in \DPD^*_{2n}$, for  each  step  $(s_{k-1},s_k)$ at height $h_k$ and $\xi_k$ ($1\leq k\leq {2n}$),  we associate the valuation
$$
v((s_{k-1},s_k), \xi_k)=
\begin{cases}
q^{\xi_k}p^{h_k-{\xi_k}}& \textrm{if $(s_{k-1},s_k)=\nearrow$,}\\
q^{\xi_k}p^{h_k-1-{\xi_k}}& \textrm{if $(s_{k-1},s_k)=\searrow$.}
\end{cases}
$$
Then
\eqref{eq:new-q-secant2} is equivalent to
\begin{align}\label{eq:FV2}
\sum_{\sigma \in \A_{2n}} p^{\thto\sigma}q^{\toht\sigma} = \sum_{(\gamma;\xi) \in \DPD^*_{2n}}\prod_{k=1}^{2n}v((s_{k-1},s_k), \xi_k)\qquad (n\geq 1).
\end{align}

We prove \eqref{eq:FV1} and \eqref{eq:FV2} by using the Fran\c{c}on-Viennot's bijection  $\Psi_{FV}:\A_{2n+1}\rightarrow \DPD_{2n}$ and $\Psi^*_{FV}:\A_{2n}\rightarrow \DPD^*_{2n}$, see \cite{FV79}.
For $k\in[n]$, define $\toht_k\sigma$, $\thot_k\sigma$, and $\thto_k\sigma$ by
$$
\begin{aligned}
\toht_k\sigma &= \#\set{l\in[n-1] : l+1<\sigma^{-1}_k \quad \text{and} \quad \sigma_l>k>\sigma_{l+1}},\\
\thot_k\sigma &= \#\set{l\in[n-1] : \sigma^{-1}_k<l \quad \text{and} \quad \sigma_l<k<\sigma_{l+1}},\\
\thto_k\sigma &= \#\set{l\in[n-1] : \sigma^{-1}_k<l \quad \text{and} \quad \sigma_l>k>\sigma_{l+1}}.
\end{aligned}
$$
For any $\sigma \in \A_{2n+1}$, the corresponding Dyck path diagramme $\Psi_{FV}(\sigma):=(s_0,\ldots, s_n; \xi_1, \ldots, \xi_n)$ is defined as follows: let $s_0=(0,0)$, and
for $k=1,\ldots, 2n$,
\begin{itemize}
\item $(s_{k-1},s_k)$ is $\nearrow$ and $\xi_k=\toht_k\sigma$ if $k$ is a valley, i.e., $\sigma^{-1}_k$ is even,
\item $(s_{k-1},s_k)$ is $\searrow$ and $\xi_k=\toht_k\sigma$ if $k$ is a peak, i.e., $\sigma^{-1}_k$ is odd.
\end{itemize}
By induction on $k$, we can show that the height $h_k$ of $(s_{k-1},s_k)$ is equal to the sum of $\toht_k\sigma$ and $\thot_k\sigma$, that is,
$$\toht_k\sigma + \thot_k\sigma = h_k.$$
Since $\toht_{2n+1}\sigma = \thot_{2n+1}\sigma = 0$,
\begin{align*}
\prod_{k=1}^{2n}v((s_{k-1},s_k), \xi_k)
= \prod_{k=1}^{2n} p^{\thot_k\sigma} q^{\toht_k\sigma} = p^{\thot\sigma} q^{\toht\sigma}.
\end{align*}
This proves \eqref{eq:FV1}.

For any $\sigma =\sigma_1\cdots\sigma_{2n} \in \A_{2n}$, let $\sigma^*=\sigma_1\cdots\sigma_{2n+1} \in \A_{2n+1}$ with $\sigma_{2n+1}=2n+1$. Then we define the bijection $\Psi^*_{FV}:\A_{2n}\rightarrow \DPD^*_{2n}$ by
$$\Psi^*_{FV}(\sigma) := \Psi_{FV}(\sigma^*).$$
If $h_k$ is the height of $(s_{k-1}, s_k)$, since
$$
\toht_k\sigma = \toht_k{\sigma^*} \quad\text{and}\quad
\thto_k\sigma = \begin{cases}
\thot_k{\sigma^*} &\text{if $\sigma^{-1}_k$ is even,}\\
\thot_k{\sigma^*} - 1 &\text{if $\sigma^{-1}_k$ is odd,}
\end{cases}$$
then
\begin{align*}
\toht_k\sigma + \thto_k\sigma =
\begin{cases}
h_k &\text{if $\sigma^{-1}_k$ is even,}\\
h_k-1 &\text{if $\sigma^{-1}_k$ is odd.}
\end{cases}
\end{align*}
It follows that
\begin{align*}
\prod_{k=1}^{2n}v((s_{k-1},s_k), \xi_k)
= \prod_{k=1}^{2n} p^{\thto_k\sigma} q^{\toht_k\sigma} = p^{\thto\sigma} q^{\toht\sigma}.
\end{align*}
This proves \eqref{eq:FV2}. We are done.
\end{proof}

\begin{rmk}
From the above continued fraction expansions it is clear that  the $(p,q)$-Euler number $E_n(p,q)$ is symmetric in $p$ and $q$.
This can also be seen  combinatorially as  follows:
For $\sigma=\sigma_1\dots\sigma_{n} \in \S_n$,  define its  reverse  to be  $\sigma^R=\sigma_n\dots\sigma_{1}$
and its complement to be
$\bar \sigma=(n+1-\sigma_1)\dots(n+1-\sigma_n).$ Clearly,
if  $\sigma \in \A_{2n}$, then $\bar\sigma^R \in \A_{2n}$ and  $(\ptoht,\pthto)\sigma=(\pthto,\ptoht)\bar\sigma^R$;
if $\sigma \in \A_{2n+1}$, then $\sigma^R \in \A_{2n+1}$ and $(\ptoht,\pthot)\sigma=(\pthot,\ptoht)\sigma^R.$
\end{rmk}

Letting $p=1$ in \eqref{eq:qtang1} and \eqref{eq:qsec1}, we recover the $q$-analogue of  Euler number $E_n$
studied by Chebikin~ \cite{Che08} and  Han et al. ~\cite{HRZ99}:
\begin{equation}\label{eq:chebikin}
E_n(q) =E_n(1,q)= \sum_{\sigma \in \A_n} q^{\toht\sigma},
\end{equation}
and, by Theorem~\ref{thm:pq}, the following continued fraction expansions,
which  were obtained by Chebikin~ \cite[Theorem 9.5]{Che08} and  Han et al. ~\cite[Eq. (42)]{HRZ99}, respectively.
\begin{corollary} We have
\begin{align}
\sum_{n=0}^{\infty} E_{2n+1}(q) t^{2n+1} &= \cfrac{t}{1-\cfrac{[1]_q[2]_qt^2}{1-\cfrac{[2]_q[3]_qt^2}{1-\cfrac{[3]_q[4]_qt^2}{\ddots}}}}
\label{eq:q-tangent}
\intertext{and}
\sum_{n=0}^{\infty} E_{2n}(q) t^{2n} &= \cfrac{1}{1-\cfrac{[1]_q^2t^2}{1-\cfrac{[2]_q^2t^2}{1-\cfrac{[3]_q^2t^2}{\ddots}}}}.
\label{eq:q-secant}
\end{align}
\end{corollary}

We can also define another $q$-analogue of tangent  number $E_{2n+1}$ and secant  number $E_{2n}$  by
\begin{align}\label{eq:new-q-euler}
E_{2n+1}^*(q) =E_{2n+1}(q^2,q)= \sum_{\sigma \in \A_n} q^{\toht\sigma+2\thot\sigma},
\intertext{and}
E_{2n}^*(q) =E_{2n}(q^2,q)= \sum_{\sigma \in \A_n} q^{\toht\sigma+2\thto\sigma}.
\end{align}
Since  $[n]_{q^2,q}=q^{n-1}[n]_q$, we derive  from Theorem~\ref{thm:pq} the corresponding generating functions.
\begin{corollary}\label{thm:qeuler}
We have
\begin{align}
\sum_{n=0}^{\infty} E_{2n+1}^*(q) t^{2n+1} &= \cfrac{t}{1-\cfrac{q[1]_q[2]_qt^2}{1-\cfrac{q^3[2]_q[3]_qt^2}{1-\cfrac{q^5[3]_q[4]_qt^2}{\ddots}}}},
\label{eq:new-q-tangent}
\intertext{and}
\sum_{n=0}^{\infty} E_{2n}^*(q) t^{2n} &= \cfrac{1}{1-\cfrac{q^0[1]_q^2t^2}{1-\cfrac{q^2[2]_q^2t^2}{1-\cfrac{q^4[3]_q^2t^2}{\ddots}}}}.
\label{eq:new-q-secant}
\end{align}
\end{corollary}

\section{Two equidistributed quintuple statistics} \label{sec:bijection}
Let $\sigma = \sigma_1\sigma_2 \cdots \sigma_n \in \S_n$.
If $1\le i \le n-1$ and $\sigma_i > \sigma_{i+1}$, the entry $\sigma_i$ (resp. $\sigma_{i+1}$) is
called a {\em descent} or {\em descent top} (resp. {\em descent bottom}) of $\sigma$. For all $i \in [n]$, the entry $\sigma_i$ is called a {\em nondescent}
or {\em nondescent top} (resp.  {\em nondescent bottom}) of $\sigma$, if it is not a descent or descent top (resp.  descent bottom) of $\sigma$.
For $\sigma\in \S_n$,  let  $\des\sigma$ be  the number of descents in $\sigma$. Then
the number of nondescents in $\sigma$ is  $\ndes\sigma := n - \des\sigma$.
For $i=1,\ldots,n$, the  entry $\sigma_i$ is called
a {\em left-to-right maximum} of $\sigma$ if $\sigma_i=\max\set{\sigma_1,\ldots,\sigma_i}$.
\begin{defn}\label{def:dz}
A nondescent $\sigma_i$  in $\sigma$ ($i=1,\ldots,n$)  is
called a  {\em foremaximum} of $\sigma$ if  $\sigma_i$ is a left-to-right maximum in $\sigma$.
Denote the number of foremaximum of $\sigma$ by $\fmax\sigma$.
A permutation $\sigma$ is called {\em coderangement} if $\fmax\sigma=0$. Let $\D^*_n$ be the subset of $\S_n$ of coderangements, that is,
$\D^*_n = \set{\sigma\in\S_n: \fmax\sigma=0}$.
\end{defn}
For example,  the  foremaximums of  $4157368 \in \S_8$ are 5 and 8, hence
 $$\fmax(4157368)=2\quad \text{and} \quad \D^*_4 = \set{2143, 3142, 3241, 4123, 4132, 4213, 4231, 4312, 4321}.$$
Recall that for $\sigma\in \S_n$  the statistic  $\MAD$ is defined by
\begin{equation}\label{eq:mad}
\MAD\sigma = \des\sigma + \toht\sigma + 2\thto\sigma.
\end{equation}
\begin{thm}\label{thm:preserving}
For $n\ge 1$, there is a bijection $\Phi$ on $\S_n$ such that
\begin{equation}\label{eq:equidis}
(\ndes, \fmax, \ptoht, \pthto, \MAD) \sigma = (\wex, \fix, \cros, \nest, \inv) \Phi(\sigma) \quad \textrm{for all $\sigma\in \S_n$}.
\end{equation}
\end{thm}

\begin{rmk}
Roselle~\cite{Ros68} first studied the statistics $(\suc, \ndes)$ on $\S_n$, where $\suc\sigma$ is the number of indices $i\in[n]$ such that $\sigma_{i+1}=\sigma_i+1$ with $\sigma_{n+1}=n+1$. (A nondescent is called a {\em rise} in his paper.)
It is trivial that $(\suc, \ndes)$ and $(\adj, \des +1)$ are equidistributed on $\S_n$ by $$(\suc, \ndes)\sigma = (\adj, \des +1)\bar\sigma,$$
where $\adj\sigma$ is the number of indices $i\in[n]$ such that $\sigma_{i+1} = \sigma_i -1$ with $\sigma_{n+1}=0$.
Burnstein~\cite{Bur10} has recently given a new proof of Eriksen's result~\cite{Eri10},  which states that $(\adj, \des +1)$ and $(\fix, \wex)$ are equidistributed on $\S_n$. Since our theorem implies that $(\fmax, \ndes)$ and $(\fix, \wex)$ are equidistributed on $\S_n$, $(\suc, \ndes)$ and $(\fmax, \ndes)$ are also equidistributed on $\S_n$. It would be interesting to make the involution on $\S_n$ keeping these statistics.
\end{rmk}

In the literature, there are several bijections on $\S_n$ transforming `linear statistics' to `cyclic statistics' on $\S_n$.
For example, the composition of the Fran\c{c}on-Viennot bijection $\Psi_{FV}:\S_n \to \L_n$ in \cite{FV79} and the inverse of the Foata-Zeilberer bijection $\Psi_{FZ}:\S_n \to \L_n$ in \cite{FZ90}
transforms the number of descents to the number of excedances.
A direct description  of this composition was given as $\Phi_{CSZ}$ in \cite{CSZ97}. It turns out that the composition of the mapping $\sigma\mapsto \sigma^{-1}$ and the bijection $\Phi_{CSZ}$ in \cite{CSZ97} satisfies the condition \eqref{eq:equidis}.
For the reader's convenience,   we give a direct description of  the mapping $\Phi$ as that  of $\Phi_{CSZ}$ in \cite{CSZ97}.
For $k\in[n]$, define  the {\em  right embracing number} $\thto_k\sigma$
(resp. {\em left embracing number} $\toht_k\sigma$) of
$k$ in $\sigma$ to be the number of integers $l\in[n-1]$ such that
$$\sigma^{-1}_k<l \text{ (resp. }l+1<\sigma^{-1}_k\text{)\quad and \quad }\sigma_l>k>\sigma_{l+1}.
$$
For example, $\thto_4(4723516)=2$ and $\toht_3(4723516)=1$.
An {\em inversion top number} (resp. {\em inversion bottom number}) of a letter $i$
in words $\pi$ is the number of occurrences of inversions of form $(i,j)$ (resp $(j,i)$) in $\pi$.
We now construct $\tau = \Phi(\sigma)$ in such a way that
$$
\thto_k\sigma = \nest_k\tau \quad \forall k=1,\ldots,n.
$$
Given a permutation $\sigma$, we first construct two biwords, $f \choose f'$ and $g \choose g'$, and then form the biword $\tau=\left({f \atop f'}~{g \atop g'}\right)$ by concatenating $f$ and $g$, and $f'$ and $g'$, respectively.
The word $f$ is defined as the subword of descent bottoms in $\sigma$, ordered increasingly, and $g$ is defined as the subword of nondescent bottoms in $\sigma$, also ordered increasingly.
The word $f'$ is the permutation on descent tops in $\sigma$ such that the inversion bottom number of each letter $a$ in $f'$ is the right embracing number of $a$ in $\sigma$. Similarly, the word $g'$ is the permutation on nondescent tops in $\sigma$ such that the inversion top number of each letter $b$ in $g'$ is the right embracing number of $b$ in $\sigma$. Rearranging the columns of $\tau'$, so that the bottom row is in increasing order, we obtain the permutation $\tau=\Phi(\sigma)$ as the top row of the rearranged biword.

\begin{ex}
Let $\sigma=4~1~ 2~ 7~ 9~6~5~8~3$, with right embracing numbers $1,0,0,2,0,1,1,0,0$.
Then
$$
{f\choose f'}
= \left(
{1 \atop 8}~
{3 \atop 4}~
{5 \atop 6}~
{6 \atop 9}
\right),
\quad
{g\choose g'}
= \left(
{2 \atop 1}~
{4 \atop 2}~
{7 \atop 7}~
{8 \atop 5}~
{9 \atop 3}
\right),
$$

$$
\tau'
= \left( {f \atop f'}~{g \atop g'} \right)
= \left(
{1 \atop 8}~
{3 \atop 4}~
{5 \atop 6}~
{6 \atop 9}~
{2 \atop 1}~
{4 \atop 2}~
{7 \atop 7}~
{8 \atop 5}~
{9 \atop 3}
\right)
\to
\left(
{2 \atop 1}~
{4 \atop 2}~
{9 \atop 3}~
{3 \atop 4}~
{8 \atop 5}~
{5 \atop 6}~
{7 \atop 7}~
{1 \atop 8}~
{6 \atop 9}
\right).
$$
and thus $\Phi(\sigma)=\tau = 2~4~9~3~8~5~7~1~6$.
\end{ex}


\begin{proof}[Proof of Theorem~\ref{thm:preserving}]
Given $\sigma\in \S_n$, let $\tau=\Phi(\sigma)$.
It is easily checked that every column ${i \choose j}$ in ${f \choose f'}$ satisfies $i<j$ and every column ${i \choose j}$ in ${g \choose g'}$ satisfies $i\ge j$.
It follows that each letter $\sigma_i$ in $g'$ corresponds to a nondescent $\sigma_i$ in $\sigma$ and also a weak excedance $\sigma_i$ in $\tau$, that is, for all $1\le i \le n$,
\begin{equation}\label{eq:ndes-wex}
\sigma_i<\sigma_{i+1} \quad (1\leq i\leq n-1)\quad\text{or}\quad i=n\iff \sigma_i \le \tau(\sigma_i).
\end{equation}
Hence,
$\ndes\,\sigma = \wex\,\tau$.
By definition, $\thto\sigma$ is the sum of the right embracing numbers in $\sigma$ and the sum of the
inversions in the  words $f'$ and $g'$.
If a pair $(i,j)$ is an inversion in $f'$ (resp. $g'$), then we have $\tau(i)<\tau(j)<j<i$ (resp. $j<i\le\tau(i)<\tau(j)$) and the
 pair $(i,j)$
 is a nesting in $\tau$. Thus
$$
\thto_k\sigma = \nest_k\tau, \quad \forall k\in [n].
$$
Similarly,
$$\toht_k\sigma = \cros_k\tau, \quad  \forall k\in [n].$$
Also, by \eqref{eq:inv} and \eqref{eq:mad},
\begin{align*}
\MAD\sigma&=n-\ndes\sigma+\toht\sigma+2\thto\sigma\\
&=n-\wex\tau+\cros\tau+2\nest\tau =\inv\tau.
\end{align*}
Finally, if $k$ is  a foremaximum in $\sigma$, then
$k$ must be  a nondescent top as well as a nondescent bottom, so $k$ is in the two subwords $g$ and $g'$.
By definition, $\toht_k\sigma=0$. Therefore,  the column ${k \choose k}$ appears in ${g \choose g'}$, i.e.,
 $\tau(k)=k$.
Conversely, if $\tau(k)=k$, then the column ${k \choose k}$ appears in ${g \choose g'}$ and $\toht_k\sigma=0$. It implies that $k$ is not a foremaximum in $\sigma$.
Hence $\fmax\sigma = \fix\tau$.
\end{proof}

\begin{ex}
We illustrate $\Phi$ on $\S_3$ with their statistics.
$$
\footnotesize
\begin{tabular}{ccccc|c||c|ccccc}
  $\ndes\sigma$ & $\fmax\sigma$ & $\toht\sigma$ & $\thto\sigma$ & $\MAD\sigma$ & $\sigma$ &  $\tau=\Phi(\sigma)$ & $\wex\tau$ & $\fix\tau$ & $\cros\tau$ & $\nest\tau$ & $\inv\tau$ \\ \hline
  3 & 3 & 0 & 0 & 0 & 123 & $(1)(2)(3)$ & 3 & 3 & 0 & 0 & 0 \\
  2 & 1 & 0 & 0 & 1 & 132 & $(1)(23)$   & 2 & 1 & 0 & 0 & 1 \\
  2 & 1 & 0 & 0 & 1 & 213 & $(12)(3)$   & 2 & 1 & 0 & 0 & 1 \\
  2 & 1 & 0 & 1 & 3 & 231 & $(13)(2)$   & 2 & 1 & 0 & 1 & 3 \\
  2 & 0 & 1 & 0 & 2 & 312 & $(123)$     & 2 & 0 & 1 & 0 & 2 \\
  1 & 0 & 0 & 0 & 2 & 321 & $(132)$     & 1 & 0 & 0 & 0 & 2 \\
\end{tabular}
$$
\end{ex}

\begin{rmk}
Corteel \cite[Section 5]{Cor05} gave a similar bijection $\Phi_{C}$ on $\S_n$ which transforms  the triple statistic $(\ndes,\ptoht,\pthto)$ to $(\wex,\nest,\cros)$, while
Steingr\'{\i}msson-Williams \cite[Section 5]{SW97} gave another one   which transforms  the triple statistic $(\ndes,\ptoht,\pthto)$ to $(n+1-\wex,\nest,\cros)$.
Other variations are given in \cite{SS96}.
The bijection $\Phi\circ\Phi_{C}^{-1}:\S_n \to \S_n$ preserves the number of weak excedances and exchanges the number of crossings and the number of nestings. A direct description of this bijection would be interesting.
\end{rmk}

\section{Counting permutations by  the  quintuple statistics }\label{sec:3}
By Theorem~\ref{thm:preserving} we have
\begin{equation}\label{eq:defA}
\begin{aligned}
A_{n}(x,y,p,q,s)&:=\sum_{\sigma\in \S_n} x^{\ndes\sigma} y^{\fmax\sigma} q^{\toht\sigma} p^{\thto\sigma} s^{\MAD\sigma} \\
&= \sum_{\sigma\in \S_n}
x^{\wex\sigma}  y^{\fix\sigma}q^{\cros\sigma} p^{\nest\sigma}s^{\inv\sigma}.
\end{aligned}
\end{equation}


Corteel~\cite{Cor05}  gave  the
continued fraction expansion of the generating function for
the triple statistic $(\wex, \cros,\nest)$ using Foata-Zeilberger's bijection. We
show here  that the same bijection can be used to compute the generating function for the quintuple statistic $(\wex, \fix, \cros,\nest,\inv)$.
For $\sigma\in S_{n}$, let
$$
\Inv\sigma=\set{(i,j):\text{$i<j$ and $\sigma_j < \sigma_i$}}.
$$
For  $k \in [n]$,
we define four types of inversions of $\sigma$ related to $k$ as follows:
\begin{itemize}
\item $\Inv_k^{(1)}\sigma=\set{(i,j)\in \Inv\sigma: \text{$j \le \sigma_j$ and $k=j$}}$,
\item $\Inv_k^{(2)}\sigma=\set{(i,j)\in \Inv\sigma: \text{$\sigma_j \le i$, $\sigma_i < i$ and $k=i$} }$,
\item $\Inv_k^{(3)}\sigma=\set{(i,j)\in \Inv\sigma: \text{$\sigma_j \le i$, $i \le \sigma_i$ and $k=i$} }$,
\item $\Inv_k^{(4)}\sigma=\set{(i,j)\in \Inv\sigma: \text{$i < \sigma_j < j$ and $k=\sigma_j$}}$,
\end{itemize}
and the set of inversions related to $k$ by
$$\Inv_k\sigma := \Inv_k^{(1)}\sigma \cup \Inv_k^{(2)}\sigma \cup \Inv_k^{(3)}\sigma \cup \Inv_k^{(4)}\sigma.$$
The {\em inversion index} on $k$ of $\sigma$ is defined by $\inv_k\sigma=\# \Inv_k\sigma$.
Clearly, we have $\inv \sigma = \sum_{k=1}^{n} \inv_k \sigma$.


\begin{thm}
\label{thm:continued_fraction}
We have
\begin{equation}\label{eq:continued_fraction}
1+\sum_{n\ge1}A_n(x,y,q,p,s) t^n
=\cfrac{1}{1-b_0t-\cfrac{a_0c_1t^2}{1-b_1t-\cfrac{a_1c_2t^2}{\ddots}}},
\end{equation}
where $a_h=x s^{2h+1} [h+1]_{q,ps}$, $b_h = xyp^h s^{2h} + (1+xq) s^h [h]_{q,ps}$, and $c_h= [h]_{q,ps}$.
\end{thm}
\begin{proof}

A {\em Laguerre history} of length $n$ is a pair  $(\gamma; \xi)$, where $\gamma$ is a Motzkin path of length $n$ and $\xi=(\xi_1,\ldots,\xi_n)$ is a sequence such that $\xi_k \in \{0, 1, \ldots, h_k\}$ (resp. $\set{-h_k, \ldots, -1, 0, 1, \ldots, h_k}$, $\{0,\ldots, h_k-1\}$) if the  step $(s_{k-1},s_{k})$ is $\nearrow$ (resp. $\longrightarrow$, $\searrow$) at height $h_k$. Let $\L_n$ be the set of Laguerre histories of length $n$.
If $(\gamma;\xi)\in \L_n$, for  each  step  $(s_{k-1},s_k)$ at height $h_k$ and $\xi_k$ ($1\leq k\leq n$),  we associate the valuation
 $$
v((s_{k-1},s_k), \xi_k)=
\begin{cases}
xq^{\xi_k}(ps)^{h_k-{\xi_k}}s^{2h_k+1}& \text{if $(s_{k-1},s_k)=\nearrow$,}\\
xy (ps)^{h_k} s^{h_k}& \text{if $(s_{k-1},s_k)=\longrightarrow$ and $\xi_k=0$,}\\
xq^{\xi_k} (ps)^{h_k-\xi_k} s^{h_k}& \text{if $(s_{k-1},s_k)=\longrightarrow$ and $\xi_k>0$,}\\
q^{-\xi_k-1} (ps)^{h_k+{\xi_k}} s^{h_k} & \text{if $(s_{k-1},s_k)=\longrightarrow$ and $\xi_k<0$,}\\
q^{\xi_k}(ps)^{h_k-1-\xi_k}& \text{if $(s_{k-1},s_k)=\searrow$,}
\end{cases}
$$
and define the valuation of the Laguerre history $(\gamma;\xi)$ by
$\prod_{k=1}^nv((s_{k-1},s_k), \xi_k)$.
 Then
\eqref {eq:continued_fraction} is equivalent to
\begin{align}\label{eq:comb-cf}
A_n(x,y,q,p,s)=\sum_{(\gamma;\xi) \in \L_n}\prod_{k=1}^nv((s_{k-1},s_k), \xi_k)\qquad (n\geq 1).
\end{align}
We prove \eqref{eq:comb-cf} by using  the Foata-Zeilberger bijection  $\Psi_{FZ}:\S_n\rightarrow \L_n$, see \cite{FZ90}.
For $k\in [n]$,
introduce  the
{\em crossing index} and {\em nesting index} on $k$ of $\sigma$ by
\begin{align*}
\cros_k\sigma&=\#\set{l: l < k \le \sigma_l <\sigma_k \quad \text{or} \quad \sigma_k <\sigma_l <k < l},\\
\nest_k\sigma&=\#\set{l: l < k \le \sigma_k <\sigma_l \quad \text{or} \quad \sigma_l <\sigma_k <k < l }.\\
\end{align*}
For  any  $\sigma \in \S_n$, the corresponding Laguerre history
$\Psi_{FZ}(\sigma):=(s_0,\ldots, s_n; \xi_1,\ldots, \xi_n)$
 is defined as follows: let $s_0=(0,0)$, and
for $k=1,\ldots, n$,
\begin{itemize}
\item $(s_{k-1},s_k)$ is $\nearrow$ and $\xi_k=\cros_k\sigma$ if $k$ is a cyclic valley, i.e., $\sigma^{-1}_k>k<\sigma_k$,
\item $(s_{k-1},s_k)$ is $\longrightarrow$ and $\xi_k=0$ if $k$ is fixed,  i.e., $\sigma_k = k$,
\item $(s_{k-1},s_k)$ is $\longrightarrow$ and $\xi_k=\cros_k\sigma$ if $k$ is a cyclic double ascent, i.e., $\sigma^{-1}_k<k<\sigma_k$,
\item $(s_{k-1},s_k)$ is $\longrightarrow$ and $\xi_k=-(\cros_k\sigma+1)$ if $k$ is a cyclic double descent, i.e., $\sigma^{-1}_k>k>\sigma_k$,
\item $(s_{k-1},s_k)$ is $\searrow$ and $\xi_k=\cros_k\sigma$ if $k$ is a cyclic peak, i.e., $\sigma^{-1}_k<k>\sigma_k$.
\end{itemize}
By induction on $k$, we can show that the height of $(s_{k-1},s_k)$ is equal to the number of indices $l$ such that $l < k \le \sigma_l$ (resp. $\sigma_l < k \le l$). Thus, if $h_k$ is the height of $(s_{k-1}, s_k)$, then
\begin{align}
\cros_k\sigma+\nest_k\sigma =
\begin{cases}
h_k &\text{if $k\le\sigma_k$,}\\
h_k-1 &\text{if $\sigma_k < k$.}
\end{cases}
\label{eq:cros_nest}
\end{align}
By definition,
\begin{align}
\nest_k \sigma &= \#\Inv_k^{(1)}\sigma + \#\Inv_k^{(2)}\sigma, \\
\#\Inv_k^{(3)}\sigma &=
\begin{cases}
h_k &\text{if $\sigma^{-1}_k \le k\le \sigma_k$}\\
h_k+1 &\text{if $k < \sigma_k$ and $k < \sigma^{-1}_k$}\\
0 &\text{if $\sigma_k < k$}
\end{cases},
\intertext{and}
\#\Inv_k^{(4)}\sigma &=
\begin{cases}
h_k &\text{if $k<\sigma^{-1}_k$ },\\
0 &\text{if $\sigma^{-1}_k\le k$}.
\end{cases}
\end{align}
So,
\begin{align}
\inv_k \sigma =
\begin{cases}
\nest_k \sigma + (2h_k+1) &\text{if $(s_{k-1},s_k)=\nearrow$}\\
\nest_k \sigma + h_k &\text{if $(s_{k-1},s_k)=\longrightarrow$}\\
\nest_k \sigma &\text{if $(s_{k-1},s_k)=\searrow$}
\end{cases}.
\label{eq:inv_nest}
\end{align}
It follows that
 \begin{align*}
 \prod_{k=1}^nv((s_{k-1},s_k), \xi_k)
 &= \prod_{k=1}^nx^{\chi(k\le \sigma_k)} y^{\chi(k=\sigma_k)} q^{\cros_k\sigma} p^{\nest_k\sigma} s^{\inv_k\sigma}\\
 &=x^{\wex\sigma} y^{\fix\sigma} q^{\cros\sigma} p^{\nest\sigma} s^{\inv \sigma}.
 \end{align*}
We are done.
\end{proof}
\begin{rmk}
Two identities \eqref{eq:cros_nest} and \eqref{eq:inv_nest} yield
\begin{equation}
\inv\sigma = n - \wex\sigma + \cros\sigma + 2 \nest\sigma.
\label{eq:inv}
\end{equation}
\end{rmk}

Substituting $p\leftarrow 1$ and $r\leftarrow 1$ in \eqref{eq:continued_fraction}, we obtain
\begin{cor}
\begin{equation}\label{eq:cf-A}
1+\sum_{n\ge1}\sum_{\sigma\in \S_n}
 x^{\wex\sigma}  y^{\fix\sigma}q^{\cros\sigma}  t^n
=\cfrac{1}{1-b_0t-\cfrac{a_0c_1t^2}{1-b_1t-\cfrac{a_1c_2t^2}{\ddots}}},
\end{equation}
where $a_h=x[h+1]_{q}$, $b_h = xy +(1+xq) [h]_{q}$, and $c_h=[h]_{q}$.
\end{cor}

Substituting $y \leftarrow  y/x$, $q\leftarrow 1$, $p\leftarrow 1$, and $s\leftarrow q$ in \eqref{eq:continued_fraction}, we obtain
\begin{cor}
\begin{equation}\label{eq:cf-SZ}
1+\sum_{n\ge1}\sum_{\sigma\in \S_n}
 x^{\exc\sigma}  y^{\fix\sigma}q^{\inv\sigma}  t^n
=\cfrac{1}{1-b_0t-\cfrac{a_0c_1t^2}{1-b_1t-\cfrac{a_1c_2t^2}{\ddots}}},
\end{equation}
where $a_h=x q^h [h+1]_{q}$, $b_h = yq^{2h} + (1+x) q^h [h]_{q}$, and $c_h= q^h [h]_{q}$.
\end{cor}

\section{The analytic proof of Theorems~\ref{thm:JV} and \ref{thm:shin-zeng}} \label{sec:analytic}
We first show how to derive Theorem~\ref{thm:JV}  from \eqref{eq:cf-A} by applying the following two contraction formulas:
\begin{align}
\cfrac{1}{1-\cfrac{c_1t}{1-\cfrac{c_2t}{\ddots}}}&=
\cfrac{1}{1-c_1t-\cfrac{c_1c_2t^2}{1-(c_2+c_3)t-\cfrac{c_3c_4t^2}{\ddots}}}\label{eq:contra1}\\
&=1+\cfrac{c_1t}{1-(c_1+c_2)t-\cfrac{c_2c_3t^2}{1-(c_3+c_4)t-\cfrac{c_4c_5t^2}{\ddots}}}\label{eq:contra2}.
\end{align}
Taking $s=-1$ and $y=1$ in \eqref{eq:cf-A},  then  $b_h = - q^h$ and $a_{h-1} c_h = -[h]_q^2$. Thus, by \eqref{eq:contra1} and \eqref{eq:contra2} with $c_{2i-1}=-[i]_q$ and $c_{2i}=[i]_q$, we get
\begin{align*}
1+\sum_{n\ge1}\sum_{\sigma\in \S_n} (-1)^{\wex\sigma} q^{\cros\sigma} t^n
&=\cfrac{1}{1+t+\cfrac{[1]_q^2t^2}{1+qt+\cfrac{[2]_q^2t^2}{1+q^2t+\cfrac{[3]_q^2t^2}{\ddots}}}}\\
&= 1- \cfrac{t}{1+\cfrac{[1]_q[2]_qt^2}{1+\cfrac{[2]_q[3]_qt^2}{1+\cfrac{[3]_q[4]_qt^2}{\ddots}}}},
\end{align*}
which is equal to $1+ \sum_{n=0}^{\infty} (-1)^{n+1} E_{2n+1}(q) t^{2n+1}$ by \eqref{eq:q-tangent}.
So \eqref{eq:jv1} is proved.

Next, taking $s=-1/q$ and $y=0$ in \eqref{eq:cf-A},  then  $b_h = 0$ and $a_{h-1} c_h = -[h]_q^2 / q$. Thus
\begin{align*}
1+\sum_{n\ge1}\sum_{\sigma\in \D_n} (-1/q)^{\exc\sigma} q^{\cros\sigma} t^n
&= \cfrac{1}{1+\cfrac{[1]_q^2t^2/q}{1+\cfrac{[2]_q^2t^2/q}{1+\cfrac{[3]_q^2t^2/q}{\ddots}}}},
\end{align*}
which is equal to $ \sum_{n=0}^{\infty} (-1/q)^{n} E_{2n}(q) t^{2n}$ by \eqref{eq:q-secant}.
This proves  \eqref{eq:jv2}.


In the same way we can  prove Theorem~\ref{thm:shin-zeng}.
Taking $s=-1/q$ and $y=1$ in \eqref{eq:cf-SZ}, then $b_h = q^{2h} + q^{2h-1}- q^{h-1}$ and $a_{h-1} c_h = -q^{2h-2}[h]_q^2$. Thus, by \eqref{eq:contra1} and \eqref{eq:contra2} with $c_{2i-1}=q^{i-1}[i]_q$ and $c_{2i}=-q^{i-1}[i]_q$, we get
\begin{align*}
1+\sum_{n\ge1}\sum_{\sigma\in \S_n} (-1/q)^{\exc\sigma} q^{\inv\sigma} t^n
&=
\cfrac{1}
{1-t+\cfrac{[1]_q^2t^2}
{1-(q^2+q-1)t+\cfrac{q^2[2]_q^2t^2}
{1-(q^4+q^3-q)t+\cfrac{q^4[3]_q^2t^2}{\ddots}}}}\\
&=
1+\cfrac{t}
{1+\cfrac{q[1]_q[2]_qt^2}
{1+\cfrac{q^3[2]_q[3]_qt^2}
{1+\cfrac{q^5[3]_q[4]_qt^2}{\ddots}}}},
\end{align*}
which is equal to $1+ \sum_{n=0}^{\infty} (-1)^{n} E_{2n+1}^*(q) t^{2n+1}$ by \eqref{eq:new-q-tangent}.
So \eqref{eq:shin-zeng1} is proved.

Next, taking $s=-1$ and $y=0$ in \eqref{eq:cf-SZ}, then  $b_h = 0$ and $a_{h-1} c_h = - q^{2h-1} [h]_q^2$. Thus
\begin{align*}
1+\sum_{n\ge1}\sum_{\sigma\in \D_n} (-1)^{\exc\sigma} q^{\inv\sigma} t^n
&= \cfrac{1}
{1+\cfrac{q [1]_q^2t^2}
{1+\cfrac{q^3 [2]_q^2t^2}
{1+\cfrac{q^5 [3]_q^2t^2}{\ddots}}}}.
\end{align*}
which is equal to $\sum_{n=0}^{\infty} (-q)^{n} E_{2n}^*(q) t^{2n}$ by \eqref{eq:new-q-secant}.
This proves \eqref{eq:shin-zeng2}.

\section{The combinatorial proof of Theorem~\ref{thm:JV}} \label{sec:combinatorial}
In view of Theorem \ref{thm:preserving}, Theorem~\ref{thm:JV} is proved if we can establish the two identities, for $n\geq 1$,
\begin{align}
\sum_{\sigma\in \S_{n}} (-1)^{\ndes\sigma} q^{\toht\sigma} &=
\begin{cases}
0&\textrm{if $n$ is even},\\
(-1)^{\frac{n+1}{2}} \sum_{\sigma\in \A_{n}} q^{ \toht\sigma}&\textrm{if $n$ is odd};
\end{cases}
\label{eq:sz1}
\intertext{and}
\sum_{\sigma\in \D^*_{n}} (-1/q)^{\ndes\sigma} q^{\toht\sigma} &=
\begin{cases}
(-1/q)^{\frac{n}{2}} \sum_{\sigma\in \A_{n}} q^{ \toht\sigma} &\textrm{if $n$ is even},\\
0 &\textrm{if $n$ is odd}.
\end{cases}
\label{eq:sz2}
\end{align}

Let $\A'_n$ (resp. $\A''_n$) be the set of alternating permutations with $\sigma_{n-1} < \sigma_n$ (resp. $\sigma_{n-1} > \sigma_n$) in $\A_n$. By definition,
$$
\A'_n =
\begin{cases}
\emptyset &\text{if $n$ is even},\\
\A_n &\text{if $n$ is odd};
\end{cases}
\quad \text{and}\quad
\A''_n =
\begin{cases}
\A_n &\text{if $n$ is even},\\
\emptyset &\text{if $n$ is odd}.
\end{cases}
$$
By \eqref{eq:chebikin} and
$$\ndes\sigma=
\begin{cases}
n+1 &\text{if $\sigma\in\A_{2n+1}$}\\
n &\text{if $\sigma\in\A_{2n}$}
\end{cases},$$
the equations in \eqref{eq:sz1} and \eqref{eq:sz2} are equivalent to two identities, for $n\ge 1$,
\begin{align}
\sum_{\sigma\in \S_{n}} (-1)^{\ndes\sigma} q^{\toht\sigma}
&= \sum_{\sigma\in \A'_{n}} (-1)^{\ndes\sigma} q^{\toht\sigma}; \label{eq:invol_phi}\\
\intertext{and}
\sum_{\sigma\in \D^*_{n}} (-1/q)^{\ndes\sigma} q^{\toht\sigma}
&= \sum_{\sigma\in \A''_{n}} (-1/q)^{\ndes\sigma} q^{\toht\sigma}.\label{eq:invol_psi}
\end{align}
We shall construct a sign-reversing involution $\phi$ on $\S_n$ (resp. $\psi$ on $\D^*_n$), which yields \eqref{eq:invol_phi} (resp. \eqref{eq:invol_psi}). Letting $\sigma=\sigma_1\cdots\sigma_n\in\S_n$, an entry $\sigma_i$ is called a {\em double ascent} (resp. {\em double descent}) if $\sigma_{i-1}<\sigma_i<\sigma_{i+1}$ (resp. $\sigma_{i-1}>\sigma_i>\sigma_{i+1}$).
\subsection{Involution $\phi$ on $\S_n$}\label{sec:phi}
Let $\sigma=\sigma_1\ldots\sigma_n \in \S_n$ with $\sigma_0=\sigma_{n+1}=0$.
Let $\sigma_m$ be the largest integer (if any) such that $\sigma_m$ is
a double ascent  or double descent in $\sigma$.
\begin{enumerate}
\item If $\sigma_m$ is a double ascent, then define
    $$
    \phi(\sigma) := \sigma_1\ldots\sigma_{m-1} \sigma_{m+1}\ldots\sigma_{k-1} \sigma_{m}\sigma_{k}\ldots\sigma_{n},
    $$
    where $k$ is a minimum index such that $\sigma_k<\sigma_m$ and $m<k$. Thus
    $$
    \toht\phi(\sigma)=\toht\sigma \quad\text{and}\quad \ndes\phi(\sigma)=\ndes\sigma-1.
    $$
\item If $\sigma_m$ is a double descent, then define
    $$
    \phi(\sigma) := \sigma_1\ldots\sigma_{k} \sigma_{m}\sigma_{k+1}\ldots\sigma_{m-1} \sigma_{m+1}\ldots\sigma_{n},
    $$
    where $k$ is a maximum index such that $\sigma_k<\sigma_m$ and $k<m$.
  Hence
    $$
    \toht\phi(\sigma)=\toht\sigma \quad\text{and}\quad \ndes\phi(\sigma)=\ndes\sigma+1.
    $$
\item Otherwise, there exists no double ascent neither double descent in $\sigma$, so
$\sigma \in \A'_n$. Define  $\A'_n$  to be the set of fixed points of $\phi$, that is, $\phi(\sigma):=\sigma$ on $\sigma \in \A'_n$.
\end{enumerate}
Therefore, when the weights of $\sigma$ is defined by $(-1)^{\ndes\sigma}q^{\toht\sigma}$, the mapping $\phi$ is a sign-reversing involution on $\S_n$ with fixed point set $\A'_n$.  Hence \eqref{eq:invol_phi}.
\subsection{Involution $\psi$ on $\D^*_n$} \label{sec:psi}
Let $\sigma=\sigma_1\ldots\sigma_n \in \D^*_n$ with $\sigma_0=0$ and $\sigma_{n+1}=n+1$.
Let $\sigma_m$ be the largest integer (if any) such that  $\sigma_m$ is a  double ascent or double descent in $\sigma$.
\begin{enumerate}
\item If $\sigma_m$ is a double ascent, then define
    $$\psi(\sigma) := \sigma_1\ldots\sigma_{k} \sigma_{m}\sigma_{k+1}\ldots\sigma_{m-1} \sigma_{m+1}\ldots\sigma_{n},$$
    where $k$ is a maximum index such that $\sigma_k>\sigma_m$ and $k<m$. Note that $k$ must  exist  because $\sigma \in \D^*_n$
    has no foremaximum.  Thus
    $$\toht\psi(\sigma)=\toht\sigma-1 \quad\text{and}\quad \ndes\psi(\sigma)=\ndes\sigma-1.$$
\item If $\sigma_m$ is a double descent, then define
    $$\psi(\sigma) := \sigma_1\ldots\sigma_{m-1} \sigma_{m+1}\ldots\sigma_{k-1} \sigma_{m}\sigma_{k}\ldots\sigma_{n},$$
    where $k$ is a minimum index such that $\sigma_k>\sigma_m$ and $m<k$. Note that $k$ could be $n+1$.
    Therefore
    $$\toht\psi(\sigma)=\toht\sigma+1 \quad\text{and}\quad \ndes\psi(\sigma)=\ndes\sigma+1.$$
\item Otherwise, there exists no double ascents or double descents in $\sigma$, so $\sigma \in \A''_n$. Let  $\A''_n$ be
the set of fixed points of $\phi$, that is, $\psi(\sigma):=\sigma$ on $\sigma \in \A''_n$.
\end{enumerate}
Therefore, when the weights of $\sigma$ is defined by $(-1/q)^{\ndes\sigma}q^{\toht\sigma}$, the mapping $\psi$ is a sign-reversing involution on $\D^*_n$ with fixed point set $\A''_n$. So we have proved \eqref{eq:invol_psi}.
\begin{rmk}
The involution $\psi$ on $\D^*_n$  keeps track of the statistic $\MAD$, that is, $\MAD\sigma = \MAD\psi(\sigma)$ for all $\sigma \in \D^*_n$. It follows that
$$\sum_{\sigma \in \D^*_n} (-1)^{\ndes\sigma} q^{\MAD\sigma}
=\sum_{\sigma \in \A''_n} (-1)^{\ndes\sigma} q^{\MAD\sigma}.$$
By \eqref{eq:defA} and $\des\sigma = \ndes\sigma=\frac{n}2$ for all $\sigma \in \A''_n$, this is a combinatorial proof of \eqref{eq:shin-zeng2}.

The involution $\phi$ on $\S_n$ does not work for \eqref{eq:shin-zeng1}. It would be interesting to find a combinatorial proof of \eqref{eq:shin-zeng1}.
\end{rmk}

\section{A parity-independent formula for $q$-Euler numbers} \label{sec:parity-indeperndent}
In 1994 Randrianarivony and Zeng~\cite{RZ94} proved the following formula for the generating function of $E_n$:
\begin{align}\label{eq:rz}
\sum_{n\geq 0}E_{n}t^n=\sum_{m\geq 0}\frac{m!t^m}{\prod_{k=0}^{\left\lfloor m/2\right\rfloor}
(1+(m-2k+1)^2t^2)}.
\end{align}
From \eqref{eq:rz} we derive immediately a parity-independent formula for $E_n$:
\begin{multline}
E_{n}=\sum_{m=0}^n\frac{1+(-1)^{n-m}}{2}\cdot\frac{m!}{4^{\floor{m/2}}} \\
\times\sum_{k=0}^{\left\lfloor m/2\right\rfloor}
\frac{(-1)^{\frac{n-m}{2}+k}(m-2k+1)^{2\left\lfloor
n/2\right\rfloor}} {k! (\left\lfloor m/2\right\rfloor-k)! (m-2k+2)_{k} (\left\lfloor
(m+1)/2\right\rfloor-k+1)_{\left\lfloor m/2\right\rfloor-k} },
\label{eq:parityofeuler}
\end{multline}
where $(a)_k := a(a+1)\cdots(a+k-1)$.

A $q$-analogue of  \eqref{eq:rz} was
given in \cite[Eq. (14), (15) and (17)]{HRZ99} :
\begin{align}\label{eq:hrz}
\sum_{n\geq 0}E_{n}(q)t^n = \sum_{m\geq 0}\frac{q^{m+1}[m]_q!t^m}{\prod_{k=0}^{\left\lfloor m/2\right\rfloor}
(q^{m-2k+1}+[m-2k+1]_q^2t^2)}.
\end{align}
Since
\begin{multline*}
\frac{t^m}{\prod_{j=0}^{\left\lfloor m/2\right\rfloor} (q^{m-2j+1}+[m-2j+1]_q^2t^2)}\\
= \sum_{k=0}^{\left\lfloor m/2\right\rfloor}
\left( \frac{\alpha_{k}}{q^{(m-2k+1)/2}-\imath[m-2k+1]_q t} + \frac{\beta_{k}}{q^{(m-2k+1)/2}+\imath[m-2k+1]_q t} \right),
\end{multline*}
where $\imath^2=-1$,  and
\begin{align*}
\alpha_{k} &= \frac{(-\imath)^m q^{(m-2k+1)m/2} [m-2k+1]_q^{2 \left\lfloor m/2 \right\rfloor-m}}
{ 2q^{(m-2k+1)/2} \prod_{j\neq k} \left( q^{m-2j+1}[m-2k+1]_q^2-q^{m-2k+1}[m-2j+1]_q^2 \right) },\\
\beta_{k} &= \frac{\imath^m q^{(m-2k+1)m/2} [m-2k+1]_q^{2\left\lfloor m/2\right\rfloor-m}}
{ 2q^{(m-2k+1)/2} \prod_{j\neq k} \left( q^{m-2j+1}[m-2k+1]_q^2-q^{m-2k+1}[m-2j+1]_q^2 \right) },
\end{align*}
we obtain then the following parity-independent formula for the $q$-Euler number:
\begin{multline*}
E_{n}(q)
=\sum_{m=0}^{n}
\frac{(-\imath)^m\imath^n+\imath^m(-\imath)^n}{2}
~q^{m+1}[m]_q! \\
\times
\sum_{k=0}^{\floor{m/2}}
\frac{q^{(m-2k+1)(m-n-2)/2} [m-2k+1]_q^{n-m+2 \floor{m/2}}}
{ \prod_{j\neq k} \left( q^{m-2j+1}[m-2k+1]_q^2-q^{m-2k+1}[m-2j+1]_q^2 \right) }.
\end{multline*}
Rewriting the denominator yields  the following result. We omit the details.
\begin{thm} \label{thm:qparityofeuler} We have
\begin{multline*}
E_n(q)
=\sum_{m=0}^n\frac{1+(-1)^{n-m}}{2}~[m]_q!\\
\times \sum_{k=0}^{\floor{m/2}}
\frac{(-1)^{\frac{n-m}{2}+k}
(1-q)^{2\floor{m/2}}
[m-2k+1]_q^{2\floor{n/2}}
q^{A_{m,k}}}
{(q^2;q^2)_k (q^2;q^2)_{\floor{m/2}-k} (q^{2(m-2k+2)};q^2)_k (q^{2(\floor{(m+1)/2}-k+1)};q^2)_{\floor{m/2}-k}},
\end{multline*}
where
\begin{align*}
A_{m,k}=k^2+k(n-m+1) - \frac{n-m}{2} + \floor{\frac{m}{2}}^2 -m \floor{\frac{n}{2}}.
\end{align*}
\end{thm}

\begin{rmk}
Clearly Theorem~\ref{thm:qparityofeuler} reduces to \eqref{eq:parityofeuler}  when $q=1$.
In \cite{JV09} a different double-sum formula for $E_n(q)$ was given. It would be interesting to show directly that the two formulas for $E_n(q)$ are equal.
\end{rmk}

\section*{Acknowledgement}
We are grateful to the anonymous referee for his careful reading and helpful remarks about a previous version of this paper.
This work is supported by la R\'egion Rh\^one-Alpes through the program ``MIRA Recherche 2008'', project 08 034147 01.


\begin{thebibliography}{SW07b}

\bibitem[And79]{And79}
D.~Andr{\'e}, \emph{D\'eveloppement de $\sec x$ et $\tan x$}, C. R. Math. Acad.
  Sci. Paris \textbf{88} (1879), 965--979.

\bibitem[Bur10]{Bur10}
A.~Burstein, \emph{On joint distribution of some permutation statistics}, 2010,
  preprint.

\bibitem[Che08]{Che08}
D.~Chebikin, \emph{Variations on descents and inversions in permutations},
  Electron. J. Combin. \textbf{15} (2008), no.~1, Research Paper 132, 34.

\bibitem[Cor05]{Cor05}
S.~Corteel, \emph{Crossings and alignments of permutations}, Adv. in Appl. Math. 38 (2007), no. 2, 149--163.

\bibitem[CSZ97]{CSZ97}
R.~J. Clarke, E.~Steingr{\'{\i}}msson, and J.~Zeng, \emph{New
  {E}uler-{M}ahonian statistics on permutations and words}, Adv. in Appl. Math.
  \textbf{18} (1997), no.~3, 237--270.

\bibitem[Eri10]{Eri10}
N.~Eriksen, \emph{Pattern and position based permutation statistics}, 2010,
  preprint.

 \bibitem[Eu55]{Eu55}  L. Euler, \emph{Institutiones calculi
 differentialis cum eius usu in analysi finitorum ac Doc- trina serierum},
  Academiae Imperialis Scientiarum Petropolitanae, St. Petersbourg, chap. VII
(Methodus summandi superior ulterius promota), 1755.


\bibitem[FH10]{FH10}
D.~Foata and G.-N. Han, \emph{The $q$-tangent and $q$-secant numbers via basic eulerian
  polynomials}, Proc. Amer. Math. Soc. \textbf{138} (2010), no.~2, 385--393.

\bibitem[Fla80]{Fla80}
P.~Flajolet, \emph{Combinatorial aspects of continued fractions}, Discrete
  Math. \textbf{32} (1980), no.~2, 125--161.

\bibitem[FS70]{FS70}
D.~Foata and M.-P. Sch{\"u}tzenberger, \emph{Th\'eorie g\'eom\'etrique des
  polyn\^omes eul\'eriens}, Lecture Notes in Mathematics, Vol. 138,
  Springer-Verlag, Berlin, 1970.

\bibitem[FV79]{FV79}
J.~Fran{\c{c}}on and G.~Viennot, \emph{Permutations selon leurs pics, creux,
  doubles mont\'ees et double descentes, nombres d'{E}uler et nombres de
  {G}enocchi}, Discrete Math. \textbf{28} (1979), no.~1, 21--35.

\bibitem[FZ90]{FZ90}
D.~Foata and D.~Zeilberger, \emph{Denert's permutation statistic is indeed
  {E}uler-{M}ahonian}, Stud. Appl. Math. \textbf{83} (1990), no.~1, 31--59.

\bibitem[HRZ99]{HRZ99}
G.-N. Han, A.~Randrianarivony, and J.~Zeng, \emph{Un autre {$q$}-analogue des
  nombres d'{E}uler}, S\'em. Lothar. Combin. \textbf{42} (1999), Art. B42e, 22
  pp. (electronic), The Andrews Festschrift (Maratea, 1998).

\bibitem[JV10]{JV09}
M.~Josuat-Verg\`{e}s, \emph{A $q$-enumeration of alternating permutations},
  European J. Combin. (2010), no.~doi:10.1016/j.ecj.2010.01.008.

\bibitem[KZ04]{KZ04}
G.~Ksavrelof and J.~Zeng, \emph{Two involutions for signed excedance numbers},
  S\'em. Lothar. Combin. \textbf{49} (2002/04), Art. B49e, 8 pp. (electronic).

\bibitem[Ros68]{Ros68}
D.~P. Roselle, \emph{Permutations by number of rises and successions}, Proc.
  Amer. Math. Soc. \textbf{19} (1968), 8--16.

\bibitem[RZ94]{RZ94}
A.~Randrianarivony and J.~Zeng, \emph{Sur une extension des nombres d'{E}uler
  et les records des permutations alternantes}, J. Combin. Theory Ser. A
  \textbf{68} (1994), no.~1, 86--99.

\bibitem[SS96]{SS96}
R.~Simion and D.~Stanton, \emph{Octabasic {L}aguerre polynomials and
  permutation statistics}, J. Comput. Appl. Math. \textbf{68} (1996), no.~1-2,
  297--329.

\bibitem[SW07a]{SWachs07}
J.~Shareshian and M.~L. Wachs, \emph{{$q$}-{E}ulerian polynomials: excedance
  number and major index}, Electron. Res. Announc. Amer. Math. Soc. \textbf{13}
  (2007), 33--45 (electronic).

\bibitem[SW07b]{SW97}
E.~Steingr{\'{\i}}msson and L.~K. Williams, \emph{Permutation tableaux and
  permutation patterns}, J. Combin. Theory Ser. A \textbf{114} (2007), no.~2,
  211--234.

\bibitem[Wil05]{Wil05}
L.~K. Williams, \emph{Enumeration of totally positive {G}rassmann cells}, Adv.
  Math. \textbf{190} (2005), no.~2, 319--342.

\end{thebibliography}

\providecommand{\bysame}{\leavevmode\hbox to3em{\hrulefill}\thinspace}
\providecommand{\href}[2]{#2}

\end{document}